\documentclass{article}
\usepackage{latexsym, amsmath, amssymb}
\usepackage{graphicx}
\usepackage{epstopdf}
\usepackage{amsmath, amsthm, amscd, amsfonts}
\usepackage{amsmath}
\usepackage{array}
\usepackage{multirow}
\usepackage{amssymb}
\usepackage{float}
\usepackage{enumerate}
\usepackage[a4paper, total={6in, 10in}]{geometry}

\newtheorem{thm}{Theorem}[section]

\theoremstyle{definition} \newtheorem{defn}[thm]{Definition}
\theoremstyle{question} 
\theoremstyle{remark} 

\makeatletter
\newcommand*{\rom}[1]{\expandafter\@slowromancap\romannumeral #1@}
\makeatother

\def\bege{\begin{equation}} \def\ende{\end{equation}}

   \def\begr{\begin{eqnarray}}
\def\endr{\end{eqnarray}} 
 
\def\bege{\begin{equation}} \def\ende{\end{equation}}
\def\begr{\begin{eqnarray}} \def\endr{\end{eqnarray}} \def\bnum{\begin{enumerate}} \def\enum{\end{enumerate}}
\begin{document}

\begin{center}
\textbf{Edge Resolvability for Circular Ladder of Heptagons}
\end{center}

\begin{center}
Jia-Bao Liu$^{1}$, Sunny Kumar Sharma$^{2}$, Vijay Kumar Bhat$^{2,a}$ and Hassan Raza$^{3}$
\end{center}
$^{1}$School of Mathematics and Physics, Anhui Jianzhu University, Hefei 230601, China.\\
$^{2}$School of Mathematics, Faculty of Sciences, Shri Mata Vaishno Devi University, Katra-$182320$, Jammu and
Kashmir, India.\\
$^{3}$Business School, University of Shanghai for Science and Technology, Shanghai 200093, China.\\
$^{1}$liujiabaoad@163.com, $^{2}$sunnysrrm94@gmail.com, $^{a}$vijaykumarbhat2000@yahoo.com, \\
$^{3}$hassan\_raza783@yahoo.com\\\\
\textbf{Abstract}
A set $\mathbb{Y}$ of elements (vertices or edges) in space is said to be a $generator$ of a metric space if each element of the space is recognized by its distances from the elements of $\mathbb{Y}$, uniquely. The generator with minimum cardinality is known as the $basis$ of the metric space, and this cardinality is the $dimension$ of the given space. In this article, we further discuss these notions with respect to a heptagonal circular ladder. We show that for a heptagonal circular ladder $\Gamma_{n}$, the edge metric dimension is three and find that it equals its metric dimension. We also introduce a new family of the convex polytope graph (denoted by $\Delta_{n}$) from a heptagonal circular ladder and find its metric dimension. Furthermore, we prove that the minimum generator (metric and edge metric) are independent for all of these families of the convex polytopes.  \\\\
\textbf{Keywords:} Metric dimension, edge metric dimension, heptagonal circular ladder, independent resolving set, convex polytope.\\\\
\textbf{MSC(2020):} 05C12, 68R01, 68R10.\\\\
\vspace{-2.0em}
\section{Introduction}
Minimum resolving Graph theory endeavors to define the conduct of real-world distance-related systems. A set $\mathbb{Y}$ of elements (vertices or edges) in space is said to be a generator of a metric space if each element of the space is determined (or recognized) by its distances from the elements of $\mathbb{Y}$, uniquely. These days there exist several types of metric generators in networks, every single one of them studied in applied and theoretical ways, corresponding to its applications or its eminence. For instance, vertex resolving sets (or simply resolving sets) \cite{ps}, independent resolving sets \cite{irs}, strong resolving sets \cite{srs}, edge resolving sets \cite{flc}, mixed resolving sets \cite{mrs}, etc. have been introduced, and contemplated respectably.\\

For the simple connected network $\Gamma$, a metric basis $\mathbb{Y}$ uniquely recognizes all the vertices of the network $\Gamma$, by mean of distance coordinates. In the recent past, some eminent researchers in Graph theory introduces that the edges of $\Gamma$ are likewise distinguished by some subset of vertices $\mathbb{Y}_{E}$ in $\Gamma$ with respect to distances to $\mathbb{Y}_{E}$ \cite{flc}. They put forward the concept of distance between an edge and a vertex of $\Gamma$, as the minimum of the distance between the given vertex and the two vertices comprising the edge, which leads them to the new theory of the metric space regarding the edges of the network $\Gamma$. The $edge$ $metric$ $dimension$ is the minimum cardinality of the set of vertices $\mathbb{Y}_{E}$ (called as the $edge$ $metric$ $basis$) in $\Gamma$ that recognizes (or resolves) every pair of different edges of $\Gamma$, in the sense that the edges of $\Gamma$ have pairwise distinctive tuples of distances to the vertices of $\mathbb{Y}_{E}$ (denoted by $edim(\Gamma)$). Next, because of recent studies of the metric and the edge metric dimension of networks, it is illustrious to talk about their applications and important investigations.\\

In seventies, Slater \cite{ps} presented the idea of a resolving set (named it, as the locating set) and its minimum cardinality for a simple connected graph, known as the metric dimension (as the location number). Melter and Harary in 1976, independently proposed a similar idea by clarifying its diverse appropriateness \cite{fr}. The investigation on this important graph-theoretic invariant is outstanding, and the number of articles has been published from both applicability and theoretic perspectives. In view of its applicability, the metric dimension fundamentally has numerous conceivable miscellaneous applications in several areas of technology, science, and social science. We also see that these notions have application in various physical situations and we discuss a few.\\\\ \clearpage
\hspace{-5.0mm}\textbf{Robotics:} Exploration of a robot can be considered in a graph-systematic network. The exploring specialist (or navigating agent) can be supposed to be a point robot, which shifts from vertex to vertex of a graph space. For this navigating specialist, there is neither a concept of visibility nor that of direction. However, it is expected that it can detect the distances to a set of milestones (or landmarks). Evidently, if the distances to an adequately large set of milestones are known by a robot, then its position in the network is determined, uniquely. Consider a navigating specialist which is exploring a space modeled by a network, and who needs to know its present location. It can impart a sign to discover how far it is from each among a set of the fixed milestones. From this, the following problem arises:\\

In a given network, what are the least number of possible landmarks required and where should they be positioned, so the distances to the landmarks decide the robot's position in the network? In the investigation of metric geometry, the accompanying concept is standard \cite{char}. A set comprising of the minimum number of landmarks that uniquely recognize the position of the robot in a network is known as a '$metric$ $basis$', and the cardinality of this metric basis is called the '$metric$ $dimension$' of the given network (usually denoted by $dim(\Gamma)$) \cite{srb}.\\\\
\textbf{Chemistry:} An elementary problem relating to chemistry is to represent mathematically the set of molecules, compounds, and atoms in a unique way, in a pharaonic structure. In this manner, the edges and vertices of a labeled network represent the bond and atom types respectively \cite{ns}.\\\\
\textbf{Other Fields:} Other important applications of metric dimension and of resolving sets can be pursued in computer science, strategies in mastermind games, mathematical disciplines, pharmacy, game theory, and signal processing where generally a moving spectator (or observer) in a network framework may be located by calculating the distance from the point to the collection of sonar stations, which have been appropriately situated in the network, for some of these aforementioned applications of metric dimension, see references in \cite{srb, ce, mt, co, sv}. \\\\
\textbf{Geometry:} A polytope in elementary geometry is a geometric object with flat sides. When polytopes are having an additional property that they are convex sets and are contained in the $n$-dimensional space $\mathbb{R}^{n}$ (Euclidean space), then they are termed as convex polytopes. Convex polytopes assume a significant job both in different branches of arithmetic and in applied zones, most quite in linear programming. The metric dimension and the edge metric dimension of a few classes of convex polytopes have been considered in \cite{sv, itm, rss, zn}. Next, we give some properties regarding the $heptagonal$ $circular$ $ladder$.\\\\
\textbf{Heptagonal Circular Ladder:} The Heptagonal circular ladder \cite{sv}, denoted by $\Gamma_{n}$, is the graph of a convex polytope with $4n$ vertices and $5n$ edges. It can be obtained from the $Pentagonal$ $Circular$ $Ladder$ $P_{n}$ \cite{rss} by introducing new vertices between the vertices of degree three in $P_{n}$. The graph of convex polytope $\Gamma_{n}$ has vertices of degrees $2$ and $3$, with cardinality $2n$ each. It comprises of $n$ $7$-sided cycles, a $2n$-sided cycle, and an $n$-sided cycle. Next, the main motivation for considering the present study is addressed. \\

The main motivation for considering the present study is now stated. Sharma and Bhat in \cite{sv}, computed the distance-based topological index called the metric dimension for three families of the graphs of convex polytopes viz., Heptagonal circular ladder $\Gamma_{n}$, and other two families are $\Sigma_{n}$ and $\Upsilon_{n}$, which are obtained from $\Gamma_{n}$ by intercalating new edges between the vertices in $\Gamma_{n}$. With regard to the metric dimension, they raised some concerns about these families of convex polytopes. Now, when we look closely at their analysis, we discover that there is only one convex polytope family left with an unknown metric dimension obtained by inserting exactly one edge into each face of the heptagon of $\Gamma_{n}$, and we denote this family by $\Delta_{n}$. The edge metric dimension of various graphs such as the Prism graph, an Antiprism graph, the Web graph, etc. has been investigated in recent history. However, the edge metric dimension for heptagonal circular ladder $\Gamma_{n}$ has not been found yet. We encourage this manuscript to fill the study gap corresponding to the metric dimension and the edge metric dimension of a heptagonal circular ladder in order to address the above-mentioned problems. In this paper, the analysis of the edge metric dimension and the metric dimension continues to expand. \\
\clearpage
The main results obtained are as follows:

\begin{itemize}
  \item For the graph $\Delta_{n}$ (which is obtained from a heptagonal circular ladder $\Gamma_{n}$), the metric dimension is three.
  \item The edge metric dimension for a heptagonal circular ladder is three.
  \item For a heptagonal circular ladder, the metric dimension and the edge metric dimension are the same.
  \item Resolving sets for a heptagonal circular ladder $\Gamma_{n}$ and all the graphs obtained from $\Gamma_{n}$ are independent.
  \item The edge resolving set for a heptagonal circular ladder is also independent.
\end{itemize}

The remainder of this paper is structured as: Section 2 introduces some basic concepts related to the metric dimension and the edge metric dimension.
Some proven heptagonal circular ladder $\Gamma_{n}$ outcomes with respect to the metric dimension are also discussed. In Section 3, by inserting new edges between the vertices of a heptagonal circular ladder $\Gamma_{n}$, we obtain a novel plane graph family, denoted by $\Delta_{n}$, and for this, we analyze its metric dimension therein. We study the edge metric dimension of a heptagonal circular ladder in Section 4 and discuss some of its properties. In Section 5, we address the problems posed for a heptagonal circular ladder about its metric dimension in the recent past. Finally, in Section 6, the conclusion and future work of the present study are discussed.

\section{Preliminaries}

We provide some fundamental descriptions of the metric dimension and the edge metric dimensions of networks in this section and analyze the recent findings of the metric dimension for a heptagonal circular ladder and its related graphs.\par

Suppose $\Gamma=\Gamma(V,E)$ be a simple (i.e., $\Gamma$ has no loops and parallel edges), finite (i.e., the cardinality of the vertex set in $\Gamma$ is finite), connected (i.e., there must exist a path between every pair of distinct vertices in $\Gamma$), and an undirected (i.e., the edges in $\Gamma$ are bidirectional) graph, with the vertex set $V$ and the edge set $E$. The distance $d_{\Gamma}(y,z)$ between two nodes $y$, $z$ in a simple connected graph $\Gamma$ is the length of the shortest $y-z$ path between the vertices $y$ and $z$, and is equal to the minimum number of edges between $y$ and $z$ in that shortest path.\\\\
\textbf{Metric Dimension:} If for any three vertices $a$, $b$, $c$ $\in V(\Gamma)$, we have $d_{\Gamma}(a,b)\neq d_{\Gamma}(a,c)$, then the vertex $a$ is said to distinguish (or resolve) the pair of vertices $b$, $c$ $(b\neq c)$ in $V(\Gamma)$. If this condition of resolvability is fulfilled by some vertices comprising a subset $\mathbb{Y}\subseteq V(\Gamma)$ i.e., every pair of different vertices in the given undirected graph $\Gamma$ is distinguished by at least one element of $\mathbb{Y}$, then $\mathbb{Y}$ is said to be a $metric$ $generator$ (or $resolving$ $set$). The $metric$ $dimension$ of the given graph $\Gamma$ is nothing, but it is just the minimum cardinality of the metric generator $\mathbb{Y}$, and is usually denoted by $\beta(\Gamma)$. The metric generator $\mathbb{Y}$ with minimum cardinality is the metric basis for $\Gamma$. For an ordered subset of vertices $\mathbb{Y}=(\varepsilon_{1}, \varepsilon_{2}, \varepsilon_{3},...,\varepsilon_{k})$, (by $\mathbb{Y}_{E}$ for edge metric dimension) the $k$-code/coordinate/representation of node $c$ in $V(\Gamma)$ is
\vspace{-7.4mm}
\begin{center}
  \begin{eqnarray*}
  \varphi(c|\mathbb{Y})&=& (d_{\Gamma}(\varepsilon_{1},c),d_{\Gamma}(\varepsilon_{2},c), d_{\Gamma}(\varepsilon_{3},c),...,d_{\Gamma}(\varepsilon_{k},c))
  \end{eqnarray*}
\end{center}

In this respect, the set $\mathbb{Y}$ is the metric generator (or resolving set) for $\Gamma$, if $\varphi(q|\mathbb{Y})\neq \varphi(p|\mathbb{Y})$, for any pair of distinct vertices $p,q \in V(\Gamma)$.\\

\begin{defn}{\textbf{Independent resolving set:}}
A set $Y$ consisting of vertices of the graph $\Gamma$, is said to be an independent resolving set for $\Gamma$, if $Y$ is both resolving set (or metric generator) and independent \cite{irs}.
\end{defn}

One can see that the metric dimension deals with the vertices of the graph by its definition, a similar concept dealing with the edges of the graph introduced by Kelenc et al. in \cite{flc}, called the edge metric dimension of the graph $\Gamma$, which uniquely identifies the edges related to a graph $\Gamma$. \\\\
\textbf{Edge Metric Dimension:} For an edge $\varepsilon=yz$ and a vertex $x$ the distance between them is defined as:
\begin{eqnarray*}
  d_{\Gamma}(x,\varepsilon)&=& min\{d_{\Gamma}(x,y),d_{\Gamma}(x,z)
\end{eqnarray*}

A subset $\mathbb{Y}_{E}$ is called an edge metric generator for $\Gamma$, if any two different edges of $\Gamma$ are distinguish by some vertex of $\mathbb{Y}_{E}$. The edge metric generator with minimum cardinality is termed as edge metric basis and that cardinality is known as the $edge$ $metric$ $dimension$ of the graph $\Gamma$, and which is denoted by $edim(\Gamma)$ or $\beta_{E}(\Gamma)$ \cite{flc}. \\

\begin{defn}{\textbf{Independent edge resolving set:}}
A set $Y_{E}$ consisting of vertices of the graph $\Gamma$, is said to be an independent edge resolving set for $\Gamma$, if $Y_{E}$ is both edge resolving set (or metric generator) and independent.
\end{defn}
Now, Khuller et al. have shown property in \cite{srb} about the metric dimension two of a connected graph $\Gamma$ and is\\\\
\textbf{Theorem 1.} {\it \cite{srb}}
{\it Let $\mathbb{A}\subseteq V(\Gamma)$ be the metric basis for the connected graph $\Gamma$ of cardinality two i.e., $|\mathbb{A}|=2$, and say $\mathbb{A}=\{\varpi, \xi\}$.  Then, the following are true:}
\begin{itemize}
  \item {\it Between the vertices $\varpi$ and $\xi$, there exists a unique shortest path $P$.}
  \item {\it The valencies of the vertices $\varpi$ and $\xi$ can never exceed $3$.}
  \item {\it The valency of any other vertex on $P$ can never exceed $5$.}
 \end{itemize}
In \cite{sv}, Sharma and Bhat considered three families of the graphs of convex polytopes (viz., $\Gamma_{n}$, $\Sigma_{n}$, and $\Upsilon_{n}$) and demonstrated that these three families have a metric dimension equal to three and constitute plane graph classes of the same constant metric dimension. They obtained the following result for a heptagonal circular ladder:\\\\
\textbf{Theorem 2.} {\it \cite{sv}}
{\it Let $\Gamma_{n}$ be the plane graph on $4n$ vertices for the positive integer $n\geq 6$. Then, $dim(\Gamma_{n})=3$ i.e., it has location number $3$.}\\\\
For the two families of plane graphs (viz., $\Sigma_{n}$ and $\Upsilon_{n}$), which are obtained by introducing additional edges in a heptagonal circular ladder $\Gamma_{n}$, and they obtained the following results\\\\
\textbf{Theorem 3.} {\it \cite{sv}}
{\it Let $\Sigma_{n}$ be the plane graph on $4n$ vertices for the positive integer $n\geq 6$. Then, $dim(\Sigma_{n})=3$ i.e., it has location number $3$.}\\\\
\textbf{Theorem 4.} {\it \cite{sv}}
{\it Let $\Upsilon_{n}$ be the plane graph on $4n$ vertices for the positive integer $n\geq 6$. Then, $dim(\Upsilon_{n})=3$ i.e., it has location number $3$.}\\\\

In the next section, we introduce a new family of the convex polytope graph denoted by $\Delta_{n}$ from a heptagonal circular ladder $\Gamma_{n}$, by introducing additional edges between the vertices of $\Gamma_{n}$ and demonstrate that the graph $\Delta_{n}$ has metric dimension three.

\section{Vertex Resolvability of the Plane Graph $\Delta_{n}$}

In this section, we begin with the discussion on the structure of the new family of the graph of convex polytope $\Delta_{n}$, obtain from a heptagonal circular ladder $\Gamma_{n}$ \cite{sv}, and study some of its basic properties, then we investigate its metric dimension.\\\\
\textbf{The Graph Of Convex Polytope $\Delta_{n}$:} The radially symmetric plane graph $\Delta_{n}$ is obtained from the graph $\Gamma_{n}$ by the insertion of new edges in the graph $\Gamma_{n}$ between the nodes $r_{t}$ and $q_{t+1}$ $(1 \leq t \leq n $, see Fig. 1). At that point, the radially symmetrical plane graph $\Delta_{n}$ comprises $n$ $5$-sided cycles, $n$ $4$-sided cycles, a $2n$-sided cycle, and an $n$-sided cycle and it consists of $6n$ edges and $4n$ number of nodes. The graph of convex polytope $\Delta_{n}$ has vertices of degrees $2$, $3$, and $4$. The set of edges and vertices of a radially symmetrical plane graph $\Delta_{n}$ is depicted separately by $E(\Delta_{n})$ and $V(\Delta_{n})$, where $V(\Delta_{n})=V(\Gamma_{n})$ and $E(\Delta_{n})=E(\Gamma_{n})\cup \{ r_{t}q_{t+1}: 1 \leq t \leq n \}$.

\begin{center}
  \begin{figure}[h!]
  \centering
  \includegraphics[width=3.5in]{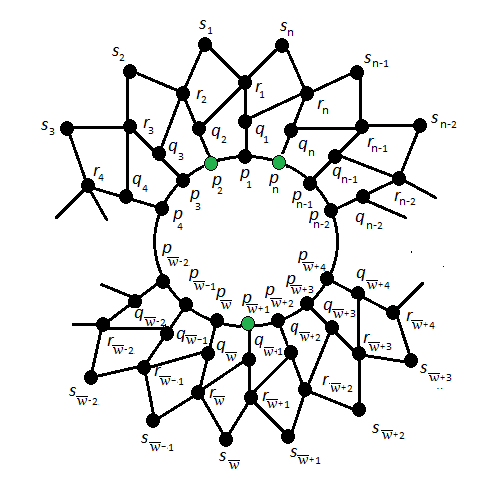}
  \caption{The Convex Polytope Graph $\Delta_{n}$}\label{p1}
\end{figure}
\end{center}

For our gentle purpose, we name the cycle generated by the set of vertices $\{p_{t}:1 \leq t \leq n\}$ in the graph, $\Delta_{n}$ as the $p$-cycle, the set of vertices $\{q_{t} : 1 \leq t \leq n\}$ in the graph, $\Delta_{n}$ as the set of mecca vertices, and the cycle generated by the set of vertices $\{r_{t}, s_{t}:1 \leq t \leq n\}$ in the graph, $\Delta_{n}$ as the $rs$-cycle. For our purpose, we consider $p_{1}=p_{n+1}$, $q_{1}=q_{n+1}$, $r_{1}=r_{n+1}$, and $s_{1}=s_{n+1}$. In the accompanying result, we determine the exact metric dimension for the family of the graph of convex polytope $\Delta_{n}$.\\\\
\textbf{Theorem 5.}
{\it Let $n$ be a positive integer such that $n \geqslant 6$ and $\Delta_{n}$ be the planar graph on $4n$ vertices as defined above. Then, we have $dim(\Delta_{n})=3$ i.e., it has location number $3$.}

\begin{proof}
In order to prove that the exact metric dimension for the family of the graph of convex polytope $\Delta_{n}$ is three, we consider two cases relying on $n$ i.e., $n \equiv 0(mod\ 2)$ and $n \equiv 1(mod\ 2)$.\\

\textbf{Case(\rom{1})} $n \equiv 0(mod\ 2)$.\\
From this, we have $n = 2\vartheta$, where $\vartheta \in \mathbb{N}$ and $\vartheta \geq 3$. Suppose $\mathbb{Y} = \{p_{2}, p_{\vartheta+1}, p_{n} \} \subset V(\Delta_{n})$ (green vertices represent the location of these three basis vertices in $\Delta_{n}$, see Fig. 1). Next, we give metric representation to every vertex of $V(\Delta_{n})\smallsetminus \mathbb{Y}$ with respect to the set $\mathbb{Y}$. \\

For, $p$-cycle $\{\nu = p_{t}:1\leq t \leq n\}$ the vertex metric representation are as follows

\begin{center}
 \begin{tabular}{|m{15.0em}|m{15.0em}|}
 \hline
 $\varphi(p_{t} |\mathbb{Y})$                                    & $\mathbb{Y}=\{p_{2}, p_{\vartheta+1}, p_{n}\}$ \\
 \hline
 $\varphi(\nu|\mathbb{Y})$: ($t=1$)                              & $(1,\vartheta,1)$ \\
 \hline
 $\varphi(\nu|\mathbb{Y})$: ($2\leq t \leq \vartheta$)           & $(t-2,\vartheta-t+1,t)$ \\
 \hline
 $\varphi(\nu|\mathbb{Y})$: ($t=\vartheta+1$)                    & $(t-2,\vartheta-t+1,2\vartheta-t)$ \\
 \hline
 $\varphi(\nu|\mathbb{Y})$: ($\vartheta+2\leq t \leq 2\vartheta$)& $(2\vartheta-t+2,t-\vartheta-1,2\vartheta-t)$ \\
 \hline
 \end{tabular}
 \end{center}

The metric representation for the set of mecca vertices $\{\nu = q_{t} : 1 \leq t \leq n\}$ are $\varphi(\nu|\mathbb{Y})=\varphi(p_{t}|\mathbb{Y})+(1,1,1)$; for $1 \leq t \leq 2\vartheta$. Finally, for $rs$-cycle $\{\nu = r_{t}, s_{t}: 1 \leq t \leq n \}$ the vertex metric representation are as follows

\begin{center}
\begin{tabular}{|m{15.0em}|m{15.0em}|}
 \hline
 $\varphi(r_{t}|\mathbb{Y})$                                     & $\mathbb{Y}=\{p_{2}, p_{\vartheta+1}, p_{n}\}$ \\
 \hline
 \vspace{-0.4mm}
 $\varphi(\nu|\mathbb{Y})$: ($t=1$)                                           & $(2,\vartheta+1,3)$ \\
 \hline
 $\varphi(\nu|\mathbb{Y})$: ($2\leq t \leq \vartheta-1$)                      & $(t,\vartheta-t+2,t+2)$ \\
 \hline
 $\varphi(\nu|\mathbb{Y})$: ($t=\vartheta$)                                   & $(\vartheta,2,\vartheta+1)$ \\
 \hline
 $\varphi(\nu|\mathbb{Y})$: ($t=\vartheta+1$)                                 & $(\vartheta+1,2,\vartheta)$ \\
 \hline

 $\varphi(\nu|\mathbb{Y})$: $(\vartheta+2\leq t \leq 2\vartheta-1)$           & $(2\vartheta-t+3,t-\vartheta+1,2\vartheta-t+1)$ \\
 \hline
 $\varphi(\nu|\mathbb{Y})$: ($t=2\vartheta$)                                  & $(2\vartheta-t+3,t-\vartheta+1,2)$ \\
 \hline
 \end{tabular}
 \end{center}
and
\begin{center}
\begin{tabular}{|m{15.0em}|m{15.0em}|}
 \hline
 $\varphi(s_{t} |\mathbb{Y})$                                     & $\mathbb{Y}=\{p_{2}, p_{\vartheta+1}, p_{n}\}$ \\
 \hline
 \vspace{-0.4mm}
 $\varphi(\nu|\mathbb{Y})$: ($t=1$)                                           & $(3,\vartheta+1,4)$ \\
 \hline
 $\varphi(\nu|\mathbb{Y})$: ($2\leq t \leq \vartheta-1$)                      & $(t+1,\vartheta-t+2,t+3)$ \\
 \hline
 $\varphi(\nu|\mathbb{Y})$: ($t=\vartheta$)                                   & $(\vartheta+1,3,\vartheta+1)$ \\
 \hline
 $\varphi(\nu|\mathbb{Y})$: ($\vartheta+1\leq t \leq 2\vartheta-2$)         & $(2\vartheta-t+3,t-\vartheta+2,2\vartheta-t+1)$ \\
 \hline
 $\varphi(\nu|\mathbb{Y})$: ($2\vartheta-1\leq t \leq 2\vartheta$)            & $(2\vartheta-t+3,t-\vartheta+2,3)$ \\
 \hline
 \end{tabular}
 \end{center}

From these codes, we see that no two elements in $V(\Delta_{n})$ have the same metric codes, suggesting $\mathbb{Y}$ to be resolving set for $\Delta_{n}$ and so, $dim(\Delta_{n})\leq 3$. Now, to complete the for this case, we have to show that $\beta(\Delta_{n}) \geq 3$ by working out that there does not exist a resolving set $\mathbb{Y}$ such that $|\mathbb{Y}|=2$. Despite what might be expected, we suppose that $\beta(\Delta_{n})= 2$. Now, by $A_{1}$, $A_{2}$, $A_{3}$, and $A_{4}$, we signify the sets with vertices as $A_{1}=\{p_{t}:1 \leq t \leq n \}$, $A_{2}=\{q_{t}:1 \leq t \leq n \}$, $A_{3}=\{r_{t} : 1 \leq t \leq n \}$, and $A_{4}=\{s_{t} : 1 \leq t \leq n\}$. At that point by Theorem $1$, we find that the degree of basis nodes can never exceed $3$. But except the vertices of the set $A_{3}$, all other vertices of the graph of the convex polytope $\Delta_{n}$ have a degree less than or equals to $3$. Therefore, we have the underlying prospects to be talked about (for $6\leq n \leq 12)$, one can check easily that no two vertices form a resolving set for $\Delta_{n}$, so we take $n\geq13$).\\\\
\begin{itemize}
\item When one, as well as the other node, is in the set $A_{l}$; $l=1,2,4.$

 \begin{center}
 \begin{tabular}{|m{13.0em}|m{18.0em}|}
 \hline
 Resolving sets & Contradictions \\
 \hline

 $\mathbb{Y}=\{p_{1},p_{g}\}$, $p_{g}$ ($2\leq g \leq n$) & $\varphi(q_{1}|\mathbb{Y})=\varphi(p_{n}|\mathbb{Y})$, for $2\leq g\leq \vartheta$ \\
 & and $\varphi(p_{2}|\mathbb{Y})=\varphi(p_{n}|\mathbb{Y})$, for $g=\vartheta+1$, a contradiction. \\
 \hline

 $\mathbb{Y}=\{q_{1},q_{g}\}$, $q_{g}$ ($2\leq g \leq n$) & $\varphi(s_{1} |\mathbb{Y}) = \varphi(s_{n} |\mathbb{Y})$, for $g=2$ \\
 & $\varphi(p_{2}|\mathbb{Y})=\varphi(q_{2}|\mathbb{Y})$, for $g=3$ \\
 & $\varphi(s_{2}|\mathbb{Y})=\varphi(r_{3}|\mathbb{Y})$, for $g=4$\\
 & $\varphi(s_{1}|\mathbb{Y})=\varphi(p_{n}|\mathbb{Y})$, for $5\leq g\leq \vartheta$\\
 & and $\varphi(p_{2}|\mathbb{Y})=\varphi(p_{n}|\mathbb{Y})$, for $g=\vartheta+1$, a contradiction. \\
 \hline

 $\mathbb{Y}=\{s_{1},s_{g}\}$, $s_{g}$ ($2\leq g \leq n$) & $\varphi(s_{n}|\mathbb{Y})=\varphi(q_{1}|\mathbb{Y})$, for $2 \leq g \leq 3$\\
 & $\varphi(s_{3}|\mathbb{Y})=\varphi(q_{4}|\mathbb{Y})$, for $4\leq g \leq 5$\\
 & and $\varphi(s_{2}|\mathbb{Y})=\varphi(q_{2}|\mathbb{Y})$, for $6\leq g\leq \vartheta+1$, a contradiction. \\
 \hline

 \end{tabular}
 \end{center}

 \item When one node is in the set $A_{1}$ and other lies in the set $A_{l}$; $l=2$ and $4.$

 \begin{center}
 \begin{tabular}{|m{13.0em}|m{18.0em}|}
 \hline
 Resolving sets & Contradictions \\
 \hline

 $\mathbb{Y}=\{p_{1}, q_{g}\}$, $q_{g}$ ($1\leq g \leq n$) & $\varphi(p_{n}|\mathbb{Y})=\varphi(p_{2}|\mathbb{Y})$, for $g=1, \vartheta+1$ \\
 & $\varphi(s_{2}|\mathbb{Y})=\varphi(s_{1}|\mathbb{Y})$, for $g=2$\\
 & and $\varphi(s_{2}|\mathbb{Y})=\varphi(q_{2}|\mathbb{Y})$, for $3 \leq g \leq \vartheta+1$, a contradiction. \\
 \hline

 $\mathbb{Y}=\{p_{1},s_{g}\}$, $s_{g}$ ($1\leq g \leq n$) & $\varphi(s_{n}|\mathbb{Y})=\varphi(q_{3}|\mathbb{Y})$, for $g=1$ \\
 & $\varphi(q_{3}|\mathbb{Y})=\varphi(s_{1}|\mathbb{Y})$, for $g=2$\\
 & and $\varphi(p_{n}|\mathbb{Y})=\varphi(q_{1}|\mathbb{Y})$, for $3 \leq g \leq \vartheta-1$\\
 & $\varphi(p_{n}|\mathbb{Y})=\varphi(p_{2}|\mathbb{Y})$, for $g=\vartheta$\\
 & and $\varphi(p_{3}|\mathbb{Y})=\varphi(q_{n}|\mathbb{Y})$, for $g=\vartheta+1$, a contradiction. \\
 \hline

 \end{tabular}
 \end{center}

 \item When one node is in the set $A_{2}$ and other lies in the set $A_{4}$.

 \begin{center}
 \begin{tabular}{|m{13.0em}|m{18.0em}|}
 \hline
 Resolving sets & Contradictions \\
 \hline

 $\mathbb{Y}=\{q_{1},s_{g}\}$, $s_{g}$ ($1\leq g \leq n$) & $\varphi(p_{n}|\mathbb{Y})=\varphi(q_{2}|\mathbb{Y})$, for $g=1$ \\
 & $\varphi(q_{2}|\mathbb{Y})=\varphi(s_{1}|\mathbb{Y})$, for $g=2$\\
 & $\varphi(s_{2}|\mathbb{Y})=\varphi(q_{3}|\mathbb{Y})$, for $3 \leq g \leq 4$\\
 & $\varphi(s_{1}|\mathbb{Y})=\varphi(p_{n}|\mathbb{Y})$, for $5 \leq g \leq \vartheta-1$\\
 & $\varphi(p_{n}|\mathbb{Y})=\varphi(p_{2}|\mathbb{Y})$, for $g=\vartheta$\\
 & and $\varphi(p_{3}|\mathbb{Y})=\varphi(r_{n-1}|\mathbb{Y})$, for $g=\vartheta+1$, a contradiction. \\
 \hline

 \end{tabular}
 \end{center}
 \end{itemize}

The above conversation reveals in this way that in this case there is no resolving set consist of two vertices for $V(\Delta_{n})$ concluding that $\beta(\Delta_{n}) = 3$.\\

\textbf{Case(\rom{2})} $n \equiv 1(mod\ 2)$.\\
From this, we have $n = 2\vartheta+1$, where $\vartheta \in \mathbb{N}$ and $\vartheta \geq 3$. Suppose $\mathbb{Y} = \{p_{2}, p_{\vartheta+1}, p_{n} \} \subset V(\Delta_{n})$. Next, we give metric representation to every vertex of $V(\Delta_{n})\smallsetminus \mathbb{Y}$ with respect to the set $\mathbb{Y}$. \\

For, $p$-cycle $\{\nu = p_{t}:1\leq t \leq n\}$ the vertex metric representation are as follows

\begin{center}
 \begin{tabular}{|m{15.0em}|m{15.0em}|}
 \hline
 $\varphi(p_{t} |\mathbb{Y})$                        & $\mathbb{Y}=\{p_{2}, p_{\vartheta+1}, p_{n}\}$ \\
 \hline
 \vspace{-0.4mm}
 $\varphi(\nu|\mathbb{Y})$:($t=1$)                                & $(1,\vartheta,1)$ \\
 \hline
 $\varphi(\nu|\mathbb{Y})$:($2\leq t \leq \vartheta$)             & $(t-2,\vartheta-t+1,t)$ \\
 \hline
 $\varphi(\nu|\mathbb{Y})$:($t=\vartheta+1$)                      & $(t-2,\vartheta-t+1,2\vartheta-t+1)$ \\
 \hline
 $\varphi(\nu|\mathbb{Y})$:($t=\vartheta+2$)                      & $(\vartheta,t-\vartheta-1,2\vartheta-t+1)$ \\
 \hline
 $\varphi(\nu|\mathbb{Y})$:($\vartheta+3\leq t \leq 2\vartheta$+1)& $(2\vartheta-t+3,t-\vartheta-1,2\vartheta-t+1)$ \\
 \hline
 \end{tabular}
 \end{center}

The metric representation for the set of mecca vertices $\{q_{t} : 1 \leq t \leq n \}$ are $\varphi(q_{t}|\mathbb{Y})=\varphi(p_{t} |\mathbb{Y})+(1,1,1)$; for $1 \leq t \leq 2\vartheta+1$. Finally, for $rs$-cycle $\{\nu=r_{t}, s_{t}: 1 \leq t \leq n\}$ the vertex metric representation are as follows

\begin{center}
 \begin{tabular}{|m{15.0em}|m{15.0em}|}
 \hline
 $\varphi(r_{t} |\mathbb{Y})$                        & $\mathbb{Y}=\{p_{2}, p_{\vartheta+1}, p_{n}\}$ \\
 \hline
 \vspace{-0.4mm}
 $\varphi(\nu|\mathbb{Y})$:($t=1$)                                & $(2,\vartheta+1,3)$ \\
 \hline
 $\varphi(\nu|\mathbb{Y})$:($2\leq t \leq \vartheta$)             & $(t,\vartheta-t+2,t+2)$ \\
 \hline
 $\varphi(\nu|\mathbb{Y})$:($t=\vartheta+1$)                      & $(\vartheta+1,2,\vartheta+1)$ \\
 \hline
 $\varphi(\nu|\mathbb{Y})$:($\vartheta+2\leq t \leq 2\vartheta$)  & $(2\vartheta-t+4,t-\vartheta+1,2\vartheta-t+2)$ \\
 \hline
 $\varphi(\nu|\mathbb{Y})$:($t= 2\vartheta$+1)                    & $(2\vartheta-t+4,t-\vartheta+1,2)$ \\
 \hline
 \end{tabular}
 \end{center}
and
\begin{center}
 \begin{tabular}{|m{15.0em}|m{15.0em}|}
 \hline
 $\varphi(s_{t} |\mathbb{Y})$                        & $\mathbb{Y}=\{p_{2}, p_{\vartheta+1}, p_{n}\}$ \\
 \hline
 \vspace{-0.4mm}
 $\varphi(\nu|\mathbb{Y})$:($t=1$)                                 & $(3,\vartheta+1,4)$ \\
 \hline
 $\varphi(\nu|\mathbb{Y})$:($2\leq t \leq \vartheta-1$)            & $(t+1,\vartheta-t+2,t+3)$ \\
 \hline
 $\varphi(\nu|\mathbb{Y})$:($t=\vartheta$)                         & $(\vartheta+1,3,\vartheta+2)$ \\
 \hline
 $\varphi(\nu|\mathbb{Y})$:($t=\vartheta+1$)                       & $(\vartheta+2,3,\vartheta+1)$ \\
 \hline
 $\varphi(\nu|\mathbb{Y})$:($\vartheta+2\leq t \leq 2\vartheta-1$) & $(2\vartheta-t+4,t-\vartheta+2,2\vartheta-t+2)$ \\
 \hline
 $\varphi(\nu|\mathbb{Y})$:($t= 2\vartheta$)                       & $(2\vartheta-t+4,t-\vartheta+2,3)$ \\
 \hline
 $\varphi(\nu|\mathbb{Y})$:($t= 2\vartheta$+1)                     & $(2\vartheta-t+4,\vartheta+2,3)$ \\
 \hline
 \end{tabular}
 \end{center}
Once again, we see that there are no two vertices with the same metric representation, suggesting that $dim(\Delta_{n}) \leq 3$. Now, in predicting that $dim(\Delta_{n}) = 2$, we consider that parallel prospects as discussed in Case(\rom{1}) are to be, and logical contradiction can be interpreted accordingly. Therefore, $dim(\Delta_{n}) = 3$ even for this case, which summarizes the theorem.
\end{proof}
Now, in view of the independent locating set, this result can also be written as:\\\\
\textbf{Theorem 6.}
{\it Let $\Delta_{n}$ be the planar graph on $4n$ vertices for the positive integer $n \geqslant 6$. Then, its independent location number is $3$.}

\begin{proof}
  For proof, refer to Theorem $5$.
\end{proof}

In the next section, we study the edge metric dimension of a heptagonal circular ladder $\Gamma_{n}$ and discuss some of its properties regarding its resolving set.

\section{Edge Resolvability of Heptagonal Circular Ladder $\Gamma_{n}$}

Some important studies were presented in \cite{flc} on the edge metric dimension of Cartesian product graphs, where the edge metric dimension value was derived for grid graphs and some cases of torus graphs (which are the Cartesian product of paths and cycles, respectively). In addition, some other findings on this aspect can be found in \cite{zn}, where this new graph parameter has been studied for some of the join graphs and for the Cartesian product of a path with any graph $\Gamma$. Furthermore, the edge metric dimension of different graphs such as the Prism graph, the Antiprism graph, the Web graphs, etc. has been studied. In relation to this fact, in recent years, other varieties of the standard metric dimensions have been deeply studied. In this context, it is now our aim to make some more contributions to this topic, and we are attempting to research the edge metric part of the circular ladder of the Heptagons. In this section, we obtain the exact edge metric dimension for a heptagonal circular ladder $\Gamma_{n}$.

\begin{center}
  \begin{figure}[h!]
  \centering
  \includegraphics[width=4in]{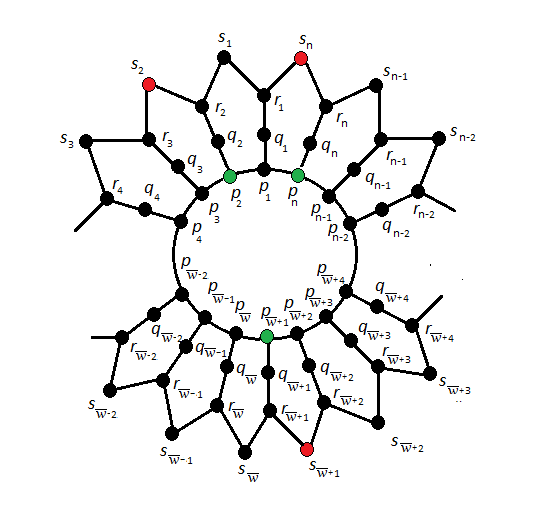}
  \caption{Heptagonal circular ladder $\Gamma_{n}$, for $n\geq6$.}\label{p1}
\end{figure}
\end{center}
\textbf{Theorem 7.}
{\it Let $\Gamma_{n}$ be a heptagonal circular ladder on $5n$ edges for the positive integer $n \geqslant 6$. Then, we have $edim(\Gamma_{n})=3$.}

\begin{proof}
It is easy to check for $6\leq n\leq 22$ that the edge metric dimension of a heptagonal circular ladder $\Gamma_{n}$ is 3 and the position  of the edge basis vertices (color in red) can be found as shown in Fig. $2$, for all $n\geq6$. Now, for $n\geq 23$, the two resulting cases depending on the positive integer $n$ are enthusiastically considered, i.e., $n \equiv 0(mod\ 2)$ and $n \equiv 1(mod\ 2)$.\\\\
\textbf{Case(\rom{1})} $n \equiv 0(mod\ 2)$.\\
From this, we have $n = 2\vartheta$, where $\vartheta \in \mathbb{N}$ and $\vartheta \geq 3$. Suppose $\mathbb{Y}_{E} = \{s_{2}, s_{\vartheta+1}, s_{n}\} \subset V(\Gamma_{n})$. Next, we give edge metric representation to every edge of $\Gamma$ with respect to the set $\mathbb{Y}_{E}$. \\\\
For, edges of $p$-cycle $\{\epsilon=p_{t}p_{t+1}| t=1,2,3,...,n\}$ the edge metric representations are as follows

\begin{center}
 \begin{tabular}{|m{15.0em}|m{15.0em}|}
 \hline
  $\varphi_{E}(\epsilon |\mathbb{Y}_{E})$ & $\mathbb{Y}_{E} = \{s_{2}, s_{\vartheta+1}, s_{n}\}$ \\
 \hline
 \vspace{-0.4mm}
 $\varphi_{E}(\epsilon)$: ($1 \leq t \leq 2$)                & $(3,\vartheta-t+3,t+2)$ \\
 \hline
 $\varphi_{E}(\epsilon)$: ($3 \leq t \leq \vartheta $)          & $(t,\vartheta-t+3,t+2)$ \\
 \hline
 $\varphi_{E}(\epsilon)$: ($t=\vartheta+1$)                     & $(t,3,2\vartheta-t+2)$ \\
 \hline
 $\varphi_{E}(\epsilon)$: ($\vartheta+2 \leq t \leq 2\vartheta-1$) & $(2\vartheta-t+4,t-\vartheta+1,2\vartheta-t+2)$ \\
 \hline
  $\varphi_{E}(\epsilon)$: ($t =2\vartheta$)                    & $(2\vartheta-t+4,t-\vartheta+1,3)$.\\
 \hline
 \end{tabular}
 \end{center}
\vspace{-0.4mm}
For the edges $\{\epsilon=p_{t}q_{t}| t=1,2,3,...,n\}$ the edge metric representations are as follows
\vspace{-0.4mm}
\begin{center}
 \begin{tabular}{|m{15.0em}|m{15.0em}|}
 \hline
  $\varphi_{E}(\epsilon |\mathbb{Y}_{E})$ & $\mathbb{Y}_{E} = \{s_{2}, s_{\vartheta+1}, s_{n}\}$ \\
 \hline
 $\varphi_{E}(\epsilon)$: ($ t=1$)                           & $(4,\vartheta+2,2)$ \\
 \hline
 $\varphi_{E}(\epsilon)$: ($2 \leq t \leq 3$)                & $(2,\vartheta-t+4,t+2)$ \\
 \hline
 $\varphi_{E}(\epsilon)$: ($4 \leq t \leq \vartheta $)          & $(t,\vartheta-t+4,t+2)$ \\
 \hline
 $\varphi_{E}(\epsilon)$: ($\vartheta+1\leq t \leq \vartheta+2$)   & $(t,2,2\vartheta-t+3)$ \\
 \hline
 $\varphi_{E}(\epsilon)$: ($\vartheta+3 \leq t \leq 2\vartheta-1$) & $(2\vartheta-t+5,t-\vartheta+1,2\vartheta-t+3)$ \\
 \hline
 \end{tabular}
 \end{center}

 \begin{center}
 \begin{tabular}{|m{15.0em}|m{15.0em}|}
 \hline
  $\varphi_{E}(\epsilon |\mathbb{Y}_{E})$ & $\mathbb{Y}_{E} = \{s_{2}, s_{\vartheta+1}, s_{n}\}$ \\
 \hline
  $\varphi_{E}(\epsilon)$: ($t =2\vartheta$)                    & $(2\vartheta-t+5,t-\vartheta+1,2)$.\\
 \hline
 \end{tabular}
 \end{center}
\vspace{-0.4mm}
For the edges $\{\epsilon=q_{t}r_{t}| t=1,2,3,...,n\}$ the edge metric representations are as follows
\vspace{-0.4mm}
\begin{center}
 \begin{tabular}{|m{7.5em}|m{11.5em}|m{7.5em}|m{11.5em}|}
 \hline
  $\varphi_{E}(\epsilon |\mathbb{Y}_{E})$ & $\mathbb{Y}_{E} = \{r_{2}, r_{\vartheta+1}, r_{n}\}$ & $\varphi_{E}(\epsilon |\mathbb{Y}_{E})$ & $\mathbb{Y}_{E} = \{r_{2}, r_{\vartheta+1}, r_{n}\}$ \\
 \hline
 \vspace{-0.4mm}
 $\varphi_{E}(\epsilon)$:($t=1$)                        & $(3,\vartheta+3,1)$   & $\varphi_{E}(\epsilon)$:($t= \vartheta+1, \vartheta+2$) &  $(t+1,1,2\vartheta-t+4)$ \\
 \hline
 $\varphi_{E}(\epsilon)$:($t=2$)  & $(1,\vartheta-t+5,3)$ & $\varphi_{E}(\epsilon)$:($t=\vartheta+3$)                  & $(n-t+6,3,n-t+4)$ \\
 \hline
 $\varphi_{E}(\epsilon)$:($t=3$)                         & $(1,\vartheta-t+5,5)$ & $\varphi_{E}(\epsilon)$:($t=\vartheta+4$)   & $(n-t+6,5,n-t+4)$ \\
 \hline
 $\varphi_{E}(\epsilon)$:($t=4$)                         & $(3,\vartheta-t+5,t+3)$ & $\varphi_{E}(\epsilon)$:($\vartheta+5 \leq t \leq n-3)$   & $(n-t+6,t-\vartheta+2, n-t+4)$ \\
 \hline
 $\varphi_{E}(\epsilon)$:($t=5$)                         & $(5,\vartheta-t+5,t+3)$ & $\varphi_{E}(\epsilon)$:($t=n-2)$   & $(n-t+6,t-\vartheta+2, 5)$ \\
 \hline
 $\varphi_{E}(\epsilon)$:($6 \leq t \leq \vartheta-2 $)     & $(t+1,\vartheta-t+5,t+3)$ & $\varphi_{E}(\epsilon)$:($t=n-1)$   & $(n-t+6,t-\vartheta+2, 3)$ \\
 \hline
 $\varphi_{E}(\epsilon)$:($t=\vartheta-1$)                  & $ (t+1,5,t+3)$ & $\varphi_{E}(\epsilon)$:($t=n)$   & $(5,t-\vartheta+2, 1)$ \\
 \hline
 $\varphi_{E}(\epsilon)$:($t=\vartheta$)                    & $ (t+1,3,t+3)$ &   & \\
 \hline
  \end{tabular}
 \end{center}
For the edges $\{\epsilon=r_{t}s_{t+1}| t=1,2,3,...,n\}$ the edge metric representations are as follows
\vspace{-0.4mm}
\begin{center}
 \begin{tabular}{|m{7.5em}|m{11.5em}|m{7.5em}|m{11.5em}|}
 \hline
  $\varphi_{E}(\epsilon |\mathbb{Y}_{E})$ & $\mathbb{Y}_{E} = \{r_{2}, r_{\vartheta+1}, r_{n}\}$ & $\varphi_{E}(\epsilon |\mathbb{Y}_{E})$ & $\mathbb{Y}_{E} = \{r_{2}, r_{\vartheta+1}, r_{n}\}$ \\
 \hline
 \vspace{-0.4mm}
 $\varphi_{E}(\epsilon)$:($t=1$)   & $(2,\vartheta+4,1)$   & $\varphi_{E}(\epsilon)$:($t=\vartheta+1$)                  & $(t+2,0,n-t+5)$ \\
 \hline
 $\varphi_{E}(\epsilon)$:($t=2$)  & $(0,\vartheta-t+6,3)$ & $\varphi_{E}(\epsilon)$:($t=\vartheta+2$)   & $(t+2,1,n-t+5)$ \\
 \hline
 $\varphi_{E}(\epsilon)$:($t=3$)  & $(1,\vartheta-t+6,5)$ & $\varphi_{E}(\epsilon)$:($t=\vartheta+3)$   & $(n-t+7,3, n-t+5)$ \\
 \hline
 $\varphi_{E}(\epsilon)$:($t=4$)  & $(3,\vartheta-t+6,7)$ & $\varphi_{E}(\epsilon)$:($t=\vartheta+4)$   & $(n-t+7,5, n-t+5)$ \\
 \hline
 $\varphi_{E}(\epsilon)$:($t=5$)  & $(5,\vartheta-t+6,t+4)$ & $\varphi_{E}(\epsilon)$:($t=\vartheta+5)$   & $(n-t+7,7, n-t+5)$ \\
 \hline
 $\varphi_{E}(\epsilon)$:($t=6$)  & $(7,\vartheta-t+6,t+4)$ &  $\varphi_{E}(\epsilon)$:($\vartheta+6 \leq t \leq n-5)$   & $(n-t+7,t-\vartheta+3, n-t+5)$ \\
 \hline
 $\varphi_{E}(\epsilon)$:($7 \leq t \leq \vartheta-4$) & $ (t+2,\vartheta-t+6,t+4)$ & $\varphi_{E}(\epsilon)$:($t=n-4)$  & $(n-t+7,t-\vartheta+3, 8)$\\
 \hline
 $\varphi_{E}(\epsilon)$:($t=\vartheta-3$)                    & $ (t+2,8,t+4)$ &$\varphi_{E}(\epsilon)$:($t=n-3)$   & $(n-t+7,t-\vartheta+3, 6)$\\
 \hline
 $\varphi_{E}(\epsilon)$:($t=\vartheta-2$)                    & $ (t+2,6,t+4)$ &$\varphi_{E}(\epsilon)$:($t=n-2)$   & $(8,t-\vartheta+3, 4)$\\
 \hline
 $\varphi_{E}(\epsilon)$:($t=\vartheta-1$)                    & $ (t+2,4,t+4)$ & $\varphi_{E}(\epsilon)$:($t=n-1)$  & $(6,t-\vartheta+3, 2)$\\
 \hline
 $\varphi_{E}(\epsilon)$:($t=\vartheta$)                      & $ (t+2,2,t+4)$ & $\varphi_{E}(\epsilon)$:($t=n)$  & $(4,t-\vartheta+3, 0)$\\
 \hline
  \end{tabular}
 \end{center}
\vspace{-0.4mm}
At last, for the edges $\{\epsilon=s_{t}r_{t+1}| t=1,2,3,...,n\}$ the edge metric representations are as follows
\vspace{-0.4mm}
\begin{center}
 \begin{tabular}{|m{7.5em}|m{11.5em}|m{7.5em}|m{11.5em}|}
 \hline
  $\varphi_{E}(\epsilon |\mathbb{Y}_{E})$ & $\mathbb{Y}_{E} = \{r_{2}, r_{\vartheta+1}, r_{n}\}$ & $\varphi_{E}(\epsilon |\mathbb{Y}_{E})$ & $\mathbb{Y}_{E} = \{r_{2}, r_{\vartheta+1}, r_{n}\}$ \\
 \hline
 \vspace{-0.4mm}
 $\varphi_{E}(\epsilon)$:($t=1$)   & $(1,\vartheta-t+5,2)$   & $\varphi_{E}(\epsilon)$:($t=\vartheta+1$)                  & $(t+3,0,n-t+4)$ \\
 \hline
 $\varphi_{E}(\epsilon)$:($t=2$)  & $(0,\vartheta-t+5,4)$ & $\varphi_{E}(\epsilon)$:($t=\vartheta+2$)   & $(n-t+6,2,n-t+4)$ \\
 \hline
 $\varphi_{E}(\epsilon)$:($t=3$)  & $(2,\vartheta-t+5,6)$ & $\varphi_{E}(\epsilon)$:($t=\vartheta+3)$   & $(n-t+6,4, n-t+4)$ \\
 \hline
 $\varphi_{E}(\epsilon)$:($t=4$)  & $(4,\vartheta-t+5,8)$ & $\varphi_{E}(\epsilon)$:($t=\vartheta+4)$   & $(n-t+6,6, n-t+4)$ \\
 \hline
 $\varphi_{E}(\epsilon)$:($t=5$)  & $(6,\vartheta-t+5,t+5)$ & $\varphi_{E}(\epsilon)$:($t=\vartheta+5)$   & $(n-t+6,8, n-t+4)$ \\
 \hline
 $\varphi_{E}(\epsilon)$:($t=6$)  & $(8,\vartheta-t+5,t+5)$ &  $\varphi_{E}(\epsilon)$:($\vartheta+6 \leq t \leq n-5)$   & $(n-t+6,t-\vartheta+4, n-t+4)$ \\
 \hline
 $\varphi_{E}(\epsilon)$:($7 \leq t \leq \vartheta-4$) & $ (t+3,\vartheta-t+5,t+5)$ & $\varphi_{E}(\epsilon)$:($t=n-4)$  & $(n-t+6,t-\vartheta+4, 7)$\\
 \hline
 $\varphi_{E}(\epsilon)$:($t=\vartheta-3$)              & $ (t+3,7,t+5)$ &$\varphi_{E}(\epsilon)$:($t=n-3)$   & $(n-t+6,t-\vartheta+4, 5)$\\
 \hline
 $\varphi_{E}(\epsilon)$:($t=\vartheta-2$)              & $ (t+3,5,t+5)$ &$\varphi_{E}(\epsilon)$:($t=n-2)$   & $(7,t-\vartheta+4, 3)$\\
 \hline
 $\varphi_{E}(\epsilon)$:($t=\vartheta-1$)              & $ (t+3,3,t+5)$ & $\varphi_{E}(\epsilon)$:($t=n-1)$  & $(5,t-\vartheta+4, 1)$\\
 \hline
 $\varphi_{E}(\epsilon)$:($t=\vartheta$)                 & $ (t+3,1,n-t+4)$ & $\varphi_{E}(\epsilon)$:($t=n)$  & $(3,t-\vartheta+4, 0)$\\
 \hline
  \end{tabular}
 \end{center}
\vspace{-0.4mm}

From these codes, we find that no two edges have the same edge metric codes, suggesting that $edim(\Gamma_{n})\leq 3$. Now, in order to finish the proof for this case, we prove that $edim(\Gamma_{n})\geq 3$ by figuring out that there is no edge metric generator $\mathbb{Y}_{E}$ such that $|\mathbb{Y}_{E}| = 2$. On the opposite, assume that $edim(\Gamma_{n})=2$. Then, we have the following prospects to be addressed at that stage.\\

\begin{center}
 \begin{tabular}{|m{13.0em}|m{23.0em}|}
 \hline
 Edge Metric Generator & Contradictions \\
 \hline

 $\mathbb{Y}_{E}=\{p_{1},p_{g}\}$, $p_{g}$ ($2\leq g \leq n$) & $\varphi_{E}(p_{1}q_{1}|\mathbb{Y}_{E})=\varphi_{E}(p_{1}p_{n}|\mathbb{Y}_{E})$, for $2\leq g \leq \vartheta$ \\
 & and $\varphi_{E}(p_{1}p_{2}|\mathbb{Y}_{E})=\varphi_{E}(p_{1}p_{n}|\mathbb{Y}_{E})$, for $g=\vartheta+1$.\\
 \hline
\end{tabular}
 \end{center}

 \begin{center}
 \begin{tabular}{|m{13.0em}|m{23.0em}|}
 \hline
 Edge Metric Generator & Contradictions \\
 \hline
 $\mathbb{Y}_{E}=\{q_{1},q_{g}\}$, $q_{g}$ ($2\leq g \leq n$) & $\varphi_{E}(p_{n}q_{n}|\mathbb{Y}_{E})=\varphi_{E}(p_{n}p_{n-1}|\mathbb{Y}_{E})$, for $2\leq g \leq \vartheta-1$ \\
 & $\varphi_{E}(r_{\vartheta-1}s_{\vartheta}|\mathbb{Y}_{E})=\varphi_{E}(r_{\vartheta}s_{\vartheta}|\mathbb{Y}_{E})$, for $g=\vartheta$\\
 & and $\varphi_{E}(p_{1}p_{2}|\mathbb{Y}_{E})=\varphi_{E}(p_{1}p_{n}|\mathbb{Y}_{E})$, for $g=\vartheta+1$.\\
 \hline

 $\mathbb{Y}_{E}=\{r_{1},r_{g}\}$, $r_{g}$ ($2\leq g \leq n$) & $\varphi_{E}(p_{n}q_{n}|\mathbb{Y}_{E})=\varphi_{E}(p_{n}p_{n-1}|\mathbb{Y}_{E})$, for $2\leq g \leq \vartheta-1$ \\
 & $\varphi_{E}(r_{\vartheta-1}s_{\vartheta}|\mathbb{Y}_{E})=\varphi_{E}(r_{\vartheta}s_{\vartheta}|\mathbb{Y}_{E})$, for $g=\vartheta$\\
 & and $\varphi_{E}(p_{1}p_{2}|\mathbb{Y}_{E})=\varphi_{E}(p_{1}p_{n}|\mathbb{Y}_{E})$, for $g=\vartheta+1$.\\
 \hline

 $\mathbb{Y}_{E}=\{s_{1},s_{g}\}$, $s_{g}$ ($2\leq g \leq n$) & $\varphi_{E}(p_{n}q_{n}|\mathbb{Y}_{E})=\varphi_{E}(p_{n}p_{n-1}|\mathbb{Y}_{E})$, for $2\leq g \leq \vartheta-1$ \\
 & $\varphi_{E}(q_{\vartheta-1}r_{\vartheta-1}|\mathbb{Y}_{E})=\varphi_{E}(p_{\vartheta}p_{\vartheta+1}|\mathbb{Y}_{E})$, for $g=\vartheta$\\
 & and $\varphi_{E}(p_{3}p_{2}|\mathbb{Y}_{E})=\varphi_{E}(p_{1}p_{n}|\mathbb{Y}_{E})$, for $g=\vartheta+1$.\\
 \hline

 $\mathbb{Y}_{E}=\{p_{1},q_{g}\}$, $q_{g}$ ($1\leq g \leq n$) & $\varphi_{E}(p_{n}q_{n}|\mathbb{Y}_{E})=\varphi_{E}(p_{n}p_{n-1}|\mathbb{Y}_{E})$, for $1\leq g \leq \vartheta-1$ \\
 & $\varphi_{E}(s_{\vartheta-1}r_{\vartheta}|\mathbb{Y}_{E})=\varphi_{E}(r_{\vartheta}s_{\vartheta}|\mathbb{Y}_{E})$, for $g=\vartheta$\\
 & and $\varphi_{E}(p_{1}p_{2}|\mathbb{Y}_{E})=\varphi_{E}(p_{1}p_{n}|\mathbb{Y}_{E})$, for $g=\vartheta+1$.\\
 \hline

 $\mathbb{Y}_{E}=\{p_{1},r_{g}\}$, $r_{g}$ ($1\leq g \leq n$) & $\varphi_{E}(p_{n}q_{n}|\mathbb{Y}_{E})=\varphi_{E}(p_{n}p_{n-1}|\mathbb{Y}_{E})$, for $1\leq g \leq \vartheta-1$ \\
 & $\varphi_{E}(s_{\vartheta-1}r_{\vartheta}|\mathbb{Y}_{E})=\varphi_{E}(r_{\vartheta}s_{\vartheta}|\mathbb{Y}_{E})$, for $g=\vartheta$\\
 & and $\varphi_{E}(p_{1}p_{2}|\mathbb{Y}_{E})=\varphi_{E}(p_{1}p_{n}|\mathbb{Y}_{E})$, for $g=\vartheta+1$.\\
 \hline

 $\mathbb{Y}_{E}=\{p_{1},s_{g}\}$, $s_{g}$ ($1\leq g \leq n$) & $\varphi_{E}(p_{n}q_{n}|\mathbb{Y}_{E})=\varphi_{E}(p_{n}p_{n-1}|\mathbb{Y}_{E})$, for $1\leq g \leq \vartheta-1$ \\
 & $\varphi_{E}(r_{\vartheta-1}q_{\vartheta-1}|\mathbb{Y}_{E})=\varphi_{E}(p_{\vartheta}p_{\vartheta+1}|\mathbb{Y}_{E})$, for $g=\vartheta$\\
 & and $\varphi_{E}(p_{\vartheta+1}p_{\vartheta+2}|\mathbb{Y}_{E})=\varphi_{E}(q_{\vartheta+3}r_{\vartheta+3}|\mathbb{Y}_{E})$, for $g=\vartheta+1$.\\
 \hline

 $\mathbb{Y}_{E}=\{q_{1},r_{g}\}$, $r_{g}$ ($1\leq g \leq n$) & $\varphi_{E}(p_{n}q_{n}|\mathbb{Y}_{E})=\varphi_{E}(p_{n}p_{n-1}|\mathbb{Y}_{E})$, for $1\leq g \leq \vartheta-1$ \\
 & $\varphi_{E}(q_{2}r_{2}|\mathbb{Y}_{E})=\varphi_{E}(p_{n-1}q_{n-1}|\mathbb{Y}_{E})$, for $g=\vartheta$\\
 & and $\varphi_{E}(p_{1}p_{2}|\mathbb{Y}_{E})=\varphi_{E}(p_{1}p_{n}|\mathbb{Y}_{E})$, for $g=\vartheta+1$.\\
 \hline

 $\mathbb{Y}_{E}=\{q_{1},s_{g}\}$, $s_{g}$ ($1\leq g \leq n$) & $\varphi_{E}(p_{n}q_{n}|\mathbb{Y}_{E})=\varphi_{E}(p_{n}p_{n-1}|\mathbb{Y}_{E})$, for $1\leq g \leq \vartheta-1$ \\
 & $\varphi_{E}(r_{\vartheta-1}q_{\vartheta-1}|\mathbb{Y}_{E})=\varphi_{E}(p_{\vartheta}p_{\vartheta+1}|\mathbb{Y}_{E})$, for $g=\vartheta$\\
 & and $\varphi_{E}(p_{\vartheta+1}p_{\vartheta+2}|\mathbb{Y}_{E})=\varphi_{E}(q_{\vartheta+3}r_{\vartheta+3}|\mathbb{Y}_{E})$, for $g=\vartheta+1$.\\
 \hline

 $\mathbb{Y}_{E}=\{r_{1},s_{g}\}$, $s_{g}$ ($1\leq g \leq n$) & $\varphi_{E}(p_{n}q_{n}|\mathbb{Y}_{E})=\varphi_{E}(p_{n}p_{n-1}|\mathbb{Y}_{E})$, for $1\leq g \leq \vartheta-1$ \\
 & $\varphi_{E}(r_{\vartheta-1}q_{\vartheta-1}|\mathbb{Y}_{E})=\varphi_{E}(p_{\vartheta}p_{\vartheta+1}|\mathbb{Y}_{E})$, for $g=\vartheta$\\
 & and $\varphi_{E}(p_{\vartheta+1}p_{\vartheta+2}|\mathbb{Y}_{E})=\varphi_{E}(q_{\vartheta+3}r_{\vartheta+3}|\mathbb{Y}_{E})$, for $g=\vartheta+1$.\\
 \hline
 \end{tabular}
 \end{center}
The above conversation reveals in this way that in this case there is no edge metric generator consist of two vertices for $\Gamma_{n}$ concluding that $edim(\Gamma_{n}) = 3$.\\\\
\textbf{Case(\rom{2})} $n \equiv 1mod\ 2)$.\\
From this, we have $n = 2\vartheta+1$, where $\vartheta \in \mathbb{N}$ and $\vartheta \geq 3$. Suppose $\mathbb{Y}_{E} = \{s_{2}, s_{\vartheta+1}, s_{n}\} \subset V(\Gamma_{n})$. Next, we give edge metric representation to every edge of $\Gamma$ with respect to the set $\mathbb{Y}_{E}$. \\\\
For, edges of $p$-cycle $\{\epsilon=p_{t}p_{t+1}| t=1,2,3,...,n\}$ the edge metric representations are as follows

\begin{center}
 \begin{tabular}{|m{15.0em}|m{15.0em}|}
 \hline
  $\varphi_{E}(\epsilon |\mathbb{Y}_{E})$ & $\mathbb{Y}_{E} = \{s_{2}, s_{\vartheta+1}, s_{n}\}$ \\
 \hline
 \vspace{-0.4mm}
 $\varphi_{E}(\epsilon)$: ($1 \leq t \leq 2$)                & $(3,\vartheta-t+3,t+2)$ \\
 \hline
 $\varphi_{E}(\epsilon)$: ($3 \leq t \leq \vartheta $)          & $(t,\vartheta-t+3,t+2)$ \\
 \hline
 $\varphi_{E}(\epsilon)$: ($t=\vartheta+1$)                     & $(t,3,2\vartheta-t+3)$ \\
 \hline
 $\varphi_{E}(\epsilon)$: ($t=\vartheta+2$)                     & $(t,t-\vartheta+1,2\vartheta-t+3)$ \\
 \hline
 $\varphi_{E}(\epsilon)$: ($\vartheta+3 \leq t \leq 2\vartheta$)   & $(2\vartheta-t+5,t-\vartheta+1,2\vartheta-t+3)$ \\
 \hline
  $\varphi_{E}(\epsilon)$: ($t =2\vartheta+1$)                  & $(2\vartheta-t+5,t-\vartheta+1,3)$.\\
 \hline
 \end{tabular}
 \end{center}
\vspace{-0.4mm}
For the edges $\{\epsilon=p_{t}q_{t}| t=1,2,3,...,n\}$, the edge metric representations are as follows
\vspace{-0.4mm}
\begin{center}
 \begin{tabular}{|m{15.0em}|m{15.0em}|}
 \hline
  $\varphi_{E}(\epsilon |\mathbb{Y}_{E})$ & $\mathbb{Y}_{E} = \{s_{2}, s_{\vartheta+1}, s_{n}\}$ \\
 \hline
 \vspace{-0.4mm}
 $\varphi_{E}(\epsilon)$: ($t=1$)                            & $(4,\vartheta-t+4,2)$ \\
 \hline
 $\varphi_{E}(\epsilon)$: ($2 \leq t \leq 3$)                & $(2,\vartheta-t+4,t+2)$ \\
 \hline
 $\varphi_{E}(\epsilon)$: ($4 \leq t \leq \vartheta $)          & $(t,\vartheta-t+4,t+2)$ \\
 \hline
 $\varphi_{E}(\epsilon)$: ($\vartheta+1\leq t \leq \vartheta+2$)   & $(t,2,2\vartheta-t+4)$ \\
 \hline
 $\varphi_{E}(\epsilon)$: ($\vartheta+3 \leq t \leq 2\vartheta$)   & $(2\vartheta-t+6,t-\vartheta+2,2\vartheta-t+4)$ \\
 \hline
  $\varphi_{E}(\epsilon)$: ($t =2\vartheta+1$)                  & $(2\vartheta-t+6,t-\vartheta+2,2)$.\\
 \hline
 \end{tabular}
 \end{center}

For the edges $\{\epsilon=q_{t}r_{t}| t=1,2,3,...,n\}$, the edge metric representations are as follows

\begin{center}
 \begin{tabular}{|m{7.5em}|m{11.5em}|m{7.5em}|m{11.5em}|}
 \hline
  $\varphi_{E}(\epsilon |\mathbb{Y}_{E})$ & $\mathbb{Y}_{E} = \{r_{2}, r_{\vartheta+1}, r_{n}\}$ & $\varphi_{E}(\epsilon |\mathbb{Y}_{E})$ & $\mathbb{Y}_{E} = \{r_{2}, r_{\vartheta+1}, r_{n}\}$ \\
 \hline
 \vspace{-0.4mm}
 $\varphi_{E}(\epsilon)$:($t=1$)                        & $(3,\vartheta-t+5,1)$   & $\varphi_{E}(\epsilon)$:($t= \vartheta+1$) &  $(t+1,1,2\vartheta-t+5)$ \\
 \hline
 $\varphi_{E}(\epsilon)$:($t=2$)  &  $(1,\vartheta-t+5,3)$ & $\varphi_{E}(\epsilon)$:($t=\vartheta+2$) &  $(t+1,1,2\vartheta-t+5)$ \\
 \hline
 $\varphi_{E}(\epsilon)$:($t=3$)                         & $(1,\vartheta-t+5,5)$ & $\varphi_{E}(\epsilon)$:($t=\vartheta+3$)                  & $(2\vartheta-t+7,3,2\vartheta-t+5)$ \\
 \hline
 $\varphi_{E}(\epsilon)$:($t=4$)                         & $(3,\vartheta-t+5,t+3)$ &$\varphi_{E}(\epsilon)$:($t=\vartheta+4$)   & $(2\vartheta-t+7,5,2\vartheta-t+5)$ \\
 \hline
 $\varphi_{E}(\epsilon)$:($t=5$)                         & $(5,\vartheta-t+5,t+3)$ &$\varphi_{E}(\epsilon)$:($\vartheta+5 \leq t \leq 2\vartheta-2)$   & $(2\vartheta-t+7,t-\vartheta+2, 2\vartheta-t+5)$ \\
 \hline
 $\varphi_{E}(\epsilon)$:($6 \leq t \leq \vartheta-2 $)     & $(t+1,\vartheta-t+5,t+3)$ &$\varphi_{E}(\epsilon)$:($t=n-2)$   & $(2\vartheta-t+7,t-\vartheta+2, 5)$ \\
 \hline
 $\varphi_{E}(\epsilon)$:($t=\vartheta-1$)                  & $ (t+1,5,t+3)$ &$\varphi_{E}(\epsilon)$:($t=n-1)$   & $(2\vartheta-t+7,t-\vartheta+2, 3)$ \\
 \hline
 $\varphi_{E}(\epsilon)$:($t=\vartheta$)                    & $ (t+1,3,t+3)$ &$\varphi_{E}(\epsilon)$:($t=n)$   & $(5,t-\vartheta+2, 1)$ \\
 \hline
  \end{tabular}
 \end{center}
For the edges $\{\epsilon=r_{t}s_{t+1}| t=1,2,3,...,n\}$, the edge metric representations are as follows

\begin{center}
 \begin{tabular}{|m{7.5em}|m{11.5em}|m{7.5em}|m{11.5em}|}
 \hline
  $\varphi_{E}(\epsilon|\mathbb{Y}_{E})$ & $\mathbb{Y}_{E} = \{r_{2}, r_{\vartheta+1}, r_{n}\}$ & $\varphi_{E}(\epsilon |\mathbb{Y}_{E})$ & $\mathbb{Y}_{E} = \{r_{2}, r_{\vartheta+1}, r_{n}\}$ \\
 \hline
 \vspace{-0.4mm}
 $\varphi_{E}(\epsilon)$:($t=1$)   & $(2,\vartheta-t+6,1)$   & $\varphi_{E}(\epsilon)$:($t=\vartheta+1$)                  & $(t+2,0,2\vartheta-t+6)$ \\
 \hline
 $\varphi_{E}(\epsilon)$:($t=2$)  & $(0,\vartheta-t+6,3)$ & $\varphi_{E}(\epsilon)$:($t=\vartheta+2$)   & $(t+2,1,2\vartheta-t+6)$ \\
 \hline
 $\varphi_{E}(\epsilon)$:($t=3$)  & $(1,\vartheta-t+6,5)$ & $\varphi_{E}(\epsilon)$:($t=\vartheta+3)$   & $(2\vartheta-t+8,3, 2\vartheta-t+6)$ \\
 \hline
 $\varphi_{E}(\epsilon)$:($t=4$)  & $(3,\vartheta-t+6,7)$ & $\varphi_{E}(\epsilon)$:($t=\vartheta+4)$   & $(2\vartheta-t+8,5, 2\vartheta-t+6)$ \\
 \hline
 $\varphi_{E}(\epsilon)$:($t=5$)  & $(5,\vartheta-t+6,t+4)$ & $\varphi_{E}(\epsilon)$:($t=\vartheta+5)$   & $(2\vartheta-t+8,7, 2\vartheta-t+6)$ \\
 \hline
 $\varphi_{E}(\epsilon)$:($t=6$)  & $(7,\vartheta-t+6,t+4)$ &  $\varphi_{E}(\epsilon)$:($\vartheta+6 \leq t \leq 2\vartheta-4)$   & $(2\vartheta-t+8,t-\vartheta+3, 2\vartheta-t+6)$ \\
 \hline
 $\varphi_{E}(\epsilon)$:($7 \leq t \leq \vartheta-4$) & $ (t+2,\vartheta-t+6,t+4)$ & $\varphi_{E}(\epsilon)$:($t=2\vartheta-3)$  & $(2\vartheta-t+8,t-\vartheta+3, 8)$\\
 \hline
 $\varphi_{E}(\epsilon)$:($t=\vartheta-3$)                    & $ (t+2,8,t+4)$ &$\varphi_{E}(\epsilon)$:($t=n-3)$   & $(2\vartheta-t+8,t-\vartheta+3, 6)$\\
 \hline
 $\varphi_{E}(\epsilon)$:($t=\vartheta-2$)                    & $ (t+2,6,t+4)$ &$\varphi_{E}(\epsilon)$:($t=n-2)$   & $(8,t-\vartheta+3, 4)$\\
 \hline
 $\varphi_{E}(\epsilon)$:($t=\vartheta-1$)                    & $ (t+2,4,t+4)$ & $\varphi_{E}(\epsilon)$:($t=n-1)$  & $(6,t-\vartheta+3, 2)$\\
 \hline
 $\varphi_{E}(\epsilon)$:($t=\vartheta$)                      & $ (t+2,2,t+4)$ & $\varphi_{E}(\epsilon)$:($t=n)$  & $(4,t-\vartheta+3, 0)$\\
 \hline
  \end{tabular}
 \end{center}

At last, for the edges $\{\epsilon=s_{t}r_{t+1}| t=1,2,3,...,n\}$, the edge metric representations are as follows

\begin{center}
 \begin{tabular}{|m{7.5em}|m{11.5em}|m{7.5em}|m{11.5em}|}
 \hline
  $\varphi_{E}(\epsilon|\mathbb{Y}_{E})$ & $\mathbb{Y}_{E} = \{r_{2}, r_{\vartheta+1}, r_{n}\}$ & $\varphi_{E}(\epsilon |\mathbb{Y}_{E})$ & $\mathbb{Y}_{E} = \{r_{2}, r_{\vartheta+1}, r_{n}\}$ \\
 \hline
 \vspace{-0.4mm}
 $\varphi_{E}(\epsilon)$:($t=1$)   & $(1,\vartheta-t+5,2)$   & $\varphi_{E}(\epsilon)$:($t=\vartheta+1$)                  & $(t+3,0,2\vartheta-t+5)$ \\
 \hline
 $\varphi_{E}(\epsilon)$:($t=2$)  & $(0,\vartheta-t+5,4)$ & $\varphi_{E}(\epsilon)$:($t=\vartheta+2$)   & $(t+3,2,2\vartheta-t+5)$ \\
 \hline
 $\varphi_{E}(\epsilon)$:($t=3$)  & $(2,\vartheta-t+5,6)$ & $\varphi_{E}(\epsilon)$:($t=\vartheta+3)$   & $(2\vartheta-t+7,4, 2\vartheta-t+5)$ \\
 \hline
 $\varphi_{E}(\epsilon)$:($t=4$)  & $(4,\vartheta-t+5,8)$ & $\varphi_{E}(\epsilon)$:($t=\vartheta+4)$   & $(2\vartheta-t+7,6, 2\vartheta-t+5)$ \\
 \hline
 $\varphi_{E}(\epsilon)$:($t=5$)  & $(6,\vartheta-t+5,t+5)$ & $\varphi_{E}(\epsilon)$:($t=\vartheta+5)$   & $(2\vartheta-t+7,8, 2\vartheta-t+5)$ \\
 \hline
 $\varphi_{E}(\epsilon)$:($t=6$)  & $(8,\vartheta-t+5,t+5)$ &  $\varphi_{E}(\epsilon)$:($\vartheta+6 \leq t \leq 2\vartheta-4)$   & $(2\vartheta-t+7,t-\vartheta+4, 2\vartheta-t+5)$ \\
 \hline
 $\varphi_{E}(\epsilon)$:($7 \leq t \leq \vartheta-4$) & $ (t+3,\vartheta-t+5,t+5)$ & $\varphi_{E}(\epsilon)$:($t=2\vartheta-3)$  & $(2\vartheta-t+7,t-\vartheta+4, 7)$\\
 \hline
 $\varphi_{E}(\epsilon)$:($t=\vartheta-3$)                    & $ (t+3,7,t+5)$ &$\varphi_{E}(\epsilon)$:($t=n-3)$   & $(2\vartheta-t+7,t-\vartheta+4, 5)$\\
 \hline
 $\varphi_{E}(\epsilon)$:($t=\vartheta-2$)                    & $ (t+3,5,t+5)$ &$\varphi_{E}(\epsilon)$:($t=n-2)$   & $(7,t-\vartheta+4, 3)$\\
 \hline
 $\varphi_{E}(\epsilon)$:($t=\vartheta-1$)                    & $ (t+3,3,t+5)$ & $\varphi_{E}(\epsilon)$:($t=n-1)$  & $(5,t-\vartheta+4, 1)$\\
 \hline
 $\varphi_{E}(\epsilon)$:($t=\vartheta$)                      & $ (t+3,1,t+5)$ & $\varphi_{E}(\epsilon)$:($t=n)$  & $(3,t-\vartheta+4, 0)$\\
 \hline
  \end{tabular}
 \end{center}

Once again, we see that there are no two vertices with the same edge metric representation, suggesting that $edim(\Gamma_{n}) \leq 3$. Now, in predicting that $edim(\Gamma_{n}) = 2$, we consider that parallel prospects as discussed in Case(\rom{1}) are to be, and logical contradiction can be interpreted accordingly. Therefore, $edim(\Gamma_{n}) = 3$ even for this case, which summarizes the theorem.
\end{proof}
For the edge metric dimension, this result can also be written as:\\\\
\textbf{Theorem 8.}
{\it Let $\Gamma_{n}$ be a heptagonal circular ladder on $5n$ edges for the positive integer $n \geqslant 6$. Then, its independent edge metric number is $3$.}

\begin{proof}
  For proof, refer to Theorem $7$.
\end{proof}

In the next part, we discuss the problems presented in \cite{sv} with the graphs corresponding to a heptagonal circular ladder with regard to its metric dimension.

\section{On Independent Resolvability of Graphs Related to Circular Ladder $\Gamma_{n}$ of Heptagons}

For two graphs viz., $\Gamma_{n}$ (i.e. heptagonal circular ladder) and $\Sigma_{n}$, Sharma and Bhat in \cite{sv} posed an open problem, that is there exist an independent resolving set for these two convex polytope graphs. In addressing these, we show in this section that there exists an independent resolving set for these two aforementioned families of the graphs of convex polytopes. Now, we have the following results for a heptagonal circular ladder:\\\\
\textbf{Theorem 9.}
{\it Let $\Gamma_{n}$ be a heptagonal circular ladder on $4n$ vertices for the positive integer $n \geqslant 6$. Then, its independent metric number is $3$.}

\begin{proof}
To show for the convex polytope graph $\Gamma_{n}$, that there exists a minimum independent resolving set $\mathbb{Y}$ of cardinality three, we follow the same technique as used in Theorem $5$.\\

Suppose $\mathbb{Y} = \{p_{2}, p_{\vartheta+1}, p_{n} \} \subset V(\Gamma_{n})$ (green vertices represent the position of these three basis vertices in $\Gamma_{n}$, see Fig. 2), $\forall$ $n\geq6$. Sharma and Bhat proved in \cite{sv} that the metric dimension for a heptagonal circular ladder is three. Now, by applying the same techniques and following the same pattern as they used (or as used in Theorem $5$), it is simple to show that the set of vertices $\mathbb{Y} = \{p_{2}, p_{\vartheta+1}, p_{n}\}$ is the independent resolving set for $\Gamma_{n}$, which concludes the theorem. \\
\end{proof}
Next, in view of the independent resolving set for the convex polytope graph $\Sigma_{n}$, we have the following result:\\\\
\textbf{Theorem 10.}
{\it Let $\Sigma_{n}$ be the convex polytope graph on $4n$ vertices for the positive integer $n \geqslant 6$. Then, its independent metric number is $3$.}

\begin{center}
  \begin{figure}[h!]
  \centering
  \includegraphics[width=3in]{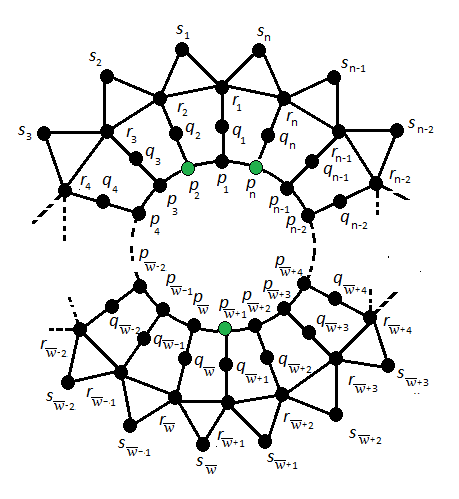}
  \caption{The Convex Polytope Graph $\Sigma_{n}$, for $n\geq6$.}\label{p1}
\end{figure}
\end{center}

\begin{proof}
To show for the convex polytope graph $\Sigma_{n}$, that there exists a minimum independent resolving set $\mathbb{Y}$ of cardinality three, we follow the same technique as used in Theorem $5$.\\

Suppose $\mathbb{Y} = \{p_{2}, p_{\vartheta+1}, p_{n} \} \subset V(\Sigma_{n})$ (green vertices represent the position of these three basis vertices in $\Sigma_{n}$, see Fig. 3), $\forall$ $n\geq6$. Sharma and Bhat proved in \cite{sv} that the metric dimension for the convex polytope graph $\Sigma_{n}$ is three. Now, by applying the same techniques and following the same pattern as they used (or as used in Theorem $5$), it is simple to show that the set of vertices $\mathbb{Y} = \{p_{2}, p_{\vartheta+1}, p_{n}\}$ is the independent resolving set for $\Sigma_{n}$, which concludes the theorem. \\
\end{proof}

Finally, all of the queries raised in Section $1$ were addressed. We conclude our work in the next section and discuss the future scope of the present study.

\section{Conclusion and Discussion}
The resolving sets for a given network or graph provide essential details required for the specific recognition of each element (vertices and edges) present in the network. In this article, we have proved that for a heptagonal circular ladder $\Gamma_{n}$, $dim(\Gamma_{n})= edim(\Gamma_{n})=3$. Thus we obtained a family of plane graphs for which the metric dimension and the edge metric dimension are the same and in this way, we have a partial answer to a question raised in \cite{flc}. We also introduced a new family of the convex polytope graph, denoted by $\Delta_{n}$, obtained from $\Gamma_{n}$ by inserting additional edges, and proved that $dim(\Delta_{n})=3$. We also addressed the problems raised for the graphs related to $\Gamma_{n}$ with regard to its metric dimension in \cite{sv}. In addition, we proved that for all of these families of convex polytopes, the minimum metric and edge metric generators are independent. The findings of this article can be useful for individuals working in the field of nano-devices, nano-engineering, and micro-devised constructed (if possible) from $\Gamma_{n}$, $\Sigma_{n}$, or $\Upsilon_{n}$. For a heptagonal circular ladder $\Gamma_{n}$ and its related graphs (viz., $\Sigma_{n}$, $\Upsilon_{n}$, and $\Delta_{n}$), we consider two distance-based graph parameters which are known as the metric dimension and the edge metric dimension, and so, regarding these, we pose some problems for future discussion:\\\\
\textbf{Problem 1:} {\it The edge metric dimension for the graphs of convex polytopes $\Sigma_{n}$, $\Upsilon_{n}$, and $\Delta_{n}$?}\\\\
\textbf{Problem 2:} {\it Whether the metric dimension and the edge metric dimension for the graphs of convex polytopes $\Sigma_{n}$, $\Upsilon_{n}$, and $\Delta_{n}$ are the same?}\\\\
\textbf{Data Availability}\\
Data sharing is not applicable to this article as no data set were generated or analyzed during the current study.\\\\
\textbf{Conflicts of Interest}\\
The authors declare no conflict of interest.\\\\
\textbf{Authors’ Contributions}\\
All the authors have equally contributed to the final manuscript.\\\\
\textbf{Acknowledgments}\\
The authors would like to thank National Natural Science Foundation of China for funding this work.

\end{document}